\documentclass[letterpaper,11pt]{article}

\usepackage[square,sort,comma,numbers]{natbib}
\usepackage{amsfonts,amsmath,verbatim,graphicx, subfigure,longtable}
\usepackage{amssymb,amsthm,amscd}
\usepackage[pdfstartview=FitH,bookmarks=false,unicode=false,pdftoolbar=false,pdfmenubar=true,]{hyperref}
\usepackage{pstricks}
\usepackage{epstopdf}
\usepackage{color}
\usepackage{transparent}
\usepackage{tikz}
\usetikzlibrary{fit, positioning}
\newtheorem{thm}{Theorem}[section]
\newtheorem{cor}[thm]{Corollary}

\newtheorem{lem}[thm]{Lemma}
\newtheorem{definition}[thm]{Definition}

\newtheorem{example}[thm]{Example}
\newtheorem{remark}[thm]{Remark}
\newcommand{\qn}{\noindent{\bf Question: \,}}
\newenvironment{pf}{\noindent\emph{Proof \,}}{\mbox{}\qed}

\numberwithin{equation}{section}
\def\R{{\mathbb R}}     % Real numbers
\def\Z{{\mathbb Z}}     % Integers
\def\P{{\mathbb P}}     % Probability
\def\E{{\mathbb E}}     % Expectation
\def\D{{\mathcal E}}    % Dirichlet form
\def\F{{\mathcal F}}    % Domain of Dirichlet form
\newcommand{\PP}{\mathbf{P}}
\newcommand{\EE}{\mathbf{E}}

\def\A{{\mathcal A}}   % Generator of a reflected diffusion
\def\eps{\varepsilon}
\def\<{{\langle}}
\def\>{{\rangle}}
\def\1{{\bf 1}}         % Indicator

\renewcommand{\bar}{\overline}

\begin{document}
\allowdisplaybreaks

\title{\Large \bf
Discrete approximations to local times \\for reflected diffusions \vspace{-0.5em}
\thanks{To appear in Electron. Commun. Probab. Vol. 21 (16), 2016. }
}
\author{{\bf Wai-Tong (Louis) Fan}}
\date{}
\maketitle

\setlength{\abovedisplayskip}{0.26cm} % Spacing above displayed equations
\setlength{\belowdisplayskip}{0.26cm} % Spacing below displayed equations

\vspace{-3em}

\begin{abstract}
We propose a discrete analogue for the boundary local time of reflected diffusions in bounded Lipschitz domains. This discrete analogue, called the discrete local time, can be effectively simulated in practice and is obtained pathwise from random walks on lattices.
%$D^{(k)}:=D\cap 2^{-k}\Z^d$ moving at rate $d\,2^{2k}$. 
We establish weak convergence of the joint law of the discrete local time and the associated random walks as the lattice size decreases to zero. 
A cornerstone of the proof is the local central limit theorem for reflected diffusions developed in \cite{zqCwtF13a}. Applications of the join convergence result to PDE problems are illustrated.
\end{abstract}

% How about Myopic scheme ?

%\medskip


%\medskip


\section{Introduction}
%\vspace{-0.5em}
Let $D\subset \R^d$ be a bounded Lipschitz domain where $d\geq 1$.  Intuitively, a \emph{reflected Brownian motion} (RBM) in $D$ is a continuous Markov process which behaves like a standard Brownian motion in the interior of $D$ and which is instantaneously pushed back by the inward normal vector $\vec{n}$ when it visits the boundary $\partial D$ of $D$.  RBMs are natural mathematical objects to study for many reasons. After all, the random motions of the pollen grains observed by Robert Brown in year 1827 were reflected at the boundary of a container. Perhaps the most notable application of RBM is the scaling limit of queuing models experiencing heavy traffic. See the monograph \cite{Harr85}. 

Reflected Brownian motion is a special case of \emph{Reflected diffusions} which we now precisely describe.
Suppose $\rho\in W^{1,2}(D)\cap\, \mathcal{C}(\bar{D})$ is a strictly positive function, and $\textbf{a}=(a^{ij})$ is a symmetric, bounded, uniformly elliptic $d\times d$ matrix-valued function with $a^{ij}\in W^{1,2}(D)$ for each $i,\,j$, where  $W^{1,2}(D):=\{f\in L^2(D):\,|\nabla f|\in L^2(D)\}$ and $\mathcal{C}(\bar{D})$ is the space of continuous functions on $\bar{D}$. It is well-known (cf.\ \cite{BH91, zqC93}) that the bilinear form $(\D,\,W^{1,2}(D))$ defined by 
$$\D(f,g) := \frac{1}{2}\int_{D}\textbf{a}\nabla f(x) \cdot \nabla g(x)\,\rho(x)\,dx$$ 
is a regular Dirichlet form in $L^2(\bar{D},\,\rho)$. Hence there is an associated Hunt process $X$ which is unique in distribution.  Furthermore, $X$ is a continuous strong Markov process in $\bar{D}$ with symmetrizing measure $\rho$ and infinitesimal generator
\begin{equation}\label{gen}
\mathcal{A} := \frac{1}{2\,\rho}\,\nabla\cdot(\rho\,\textbf{a}\nabla). 
\end{equation}

\begin{definition}\rm\label{Def:ReflectedDiffusion} The process $X$ constructed above is called an {\em $\A$-reflected diffusion}.  An important case is when $\textbf{a}$ is the identity matrix, then $X$ is called a reflected Brownian motion with drift $\frac{1}{2}\,\nabla (\log \rho)$. If in addition $\rho=1$, then $X$ is called a {\em reflected Brownian motion {\rm (}RBM\,{\rm)}}.
\end{definition}
Intuitively, $X$ behaves like a diffusion process associated to the elliptic operator $\A$ in the interior of  $D$, and  % (which is a standard Brownian motion when $\A=\frac{\Delta}{2}$);
it is instantaneously pushed back in the direction of the co-normal $\vec{\nu}:=\textbf{a}\vec{n}$ when $X_t\in \partial D$, where $\vec{n}$ is the inward unit normal.
Let $\sigma$ be  the surface measure on $\partial D$. It is well-known that there is a unique positive continuous additive functional (PCAF) of $X$ corresponding to $\sigma/2$. See, for instance, the appendix of \cite{CF12}. This PCAF  $L=(L_t)_{t\geq 0}$ is called the \emph{boundary local time} of $X$. It describes the amount of time $X$ spends near the boundary in the sense that
\begin{equation}\label{E:ConvLocal}
\lim_{\delta\to 0}\frac{1}{2\delta}\int_0^t \1\{X_s\in D^{\delta}\}\,ds = L_t \quad\text{in probability},
\end{equation}
where $D^{\delta}:= \{x\in D:\,dist(x,\,\partial D)<\delta\}$ and $\1$ is the indicator function. Moreover, $X$ admits
%, in the sense of weak solutions\footnote{Strong uniqueness for \eqref{E:SkorokhodRep} fails for some Lipschitz domains, see \cite{BBC05} for a counter example.} to stochastic differential equations, 
the Skorohod decomposition
\begin{equation}\label{E:SkorokhodRep}
X_t = X_0+ \int_0^t \beta(X_s)\cdot dB_s + \int_0^t \vec{b}(X_s)ds + \int_0^t \rho\,\textbf{a}\,\vec{n}\,(X_s)\,dL_s,\quad t\geq 0,
\end{equation}
where $B$ is the standard $d$-dimensional Brownian motion, $\vec{b}= \frac{1}{2}\,(\nabla \cdot \textbf{a}+\textbf{a}\nabla\log \rho)$ is the drift and $\beta^2 = \textbf{a}$. See \cite{BH91, zqC93, vgP90} for well-known properties about $X$ and $L$.

Discrete approximations for reflected diffusions are ubiquitous in scientific literature. However, many of them are adhoc and not rigorously justified. 
For numerical approximation schemes with mathematical justifications, there is a concise survey in  \cite[Section 5.6]{MT04}. For random walk approximation schemes, Burdzy and Chen \cite{BC08, BC11} considered  RBM for a large class of domains $D$ which contains bounded Lipschitz domains. Using Dirichlet form method and some probabilistic tools, they showed that the laws of both discrete time and continuous time simple random walks (SRW) on $D^{(k)}:= D\cap 2^{-k}\Z^d$  moving at rate $d\,2^{2k}$ converge weakly, as $k\to\infty$, to the law of RBM in $D$. The following natural question is the motivation of this paper.

\medskip

\qn{What is a discrete analogue to the boundary local time of a reflected diffusion? } 

\medskip

We consider this question interesting in its own right and in applications. 
%as important as the study of local time of reflected diffusions. Furthermore, 
A suitable candidate for such a discrete analogue, henceforth called \emph{discrete local time}, is useful in the study of partially reflected diffusions \cite{Singer08} and in stochastic particle systems in domains (such as \cite{zqCwtF13a, zqCwtF14b}) in which non-trivial interactions among particles occur only at the boundary. It can also be used to generate Monte Carlo approximations to boundary value problems in partial differential equations; see the application immediately after the statement of Theorem \ref{T:Main}. 

A rigorous answer to the above question \emph{does not} follow directly from \cite{BC08, BC11} or other published results; extra work is required to construct such an analogue and to prove convergence. To see this, we consider the case when $X$ is a RBM. Results in \cite{BC08} imply that for \textit{fixed} $\delta>0$,
\begin{equation}\label{E:Candidate1}
A^{(k)}_{\delta}(t):= \frac{1}{2\delta}\int_0^t \1\{X^{(k)}_s\in D^{\delta}\}\,ds \,\to\, \frac{1}{2\delta}\int_0^t \1\{X_s\in D^{\delta}\}\,ds
\end{equation}
in distribution as $k\to\infty$, where we used the same notation $X^{(k)}$ to denote both discrete time and continuous time SRW on $D^{(k)}$ moving at rate $d\,2^{2k}$. For discrete time SRW, the time parameter is extended by interpolation. Even though we have \eqref{E:ConvLocal}, the results in \cite{BC08} \emph{do not} tell us how small $\delta$ should be taken relative to $k$. 

One might guess that a possible candidate is the left-hand side of \eqref{E:Candidate1} with $\delta= C\,2^{-k}$ for some constant $C>0$ large enough so that for all $k\in \mathbb{N}$, we have $D^{C\,2^{-k}}$ contains the \emph{graph-boundary} $\partial D^{(k)}:= \{x\in D^{(k)}:\,v_{k}(x)<2d\}$,  where $v_{k}(x)$ is the degree of the vertex $x$ in $D^{(k)}$. 
%Such a constant $C$ can be chosen to depend only on the Lipschitz constant of $\partial D$. 
However, this candidate $A^{(k)}_{C\,2^{-k}}(t)$ turns out to be {\it problematic} since it is too sensitive to the local configuration of the graph $D^{(k)}$ near the boundary. Our result also indicates that the ``naive" candidate 
$$\; \frac{1}{2\,(2^{-k})}\int_0^t \1\{X^{(k)}_s\in \partial D^{(k)}\}\,ds,$$
which records the amount of time the random walk spends on $\partial D^{(k)}$, \emph{does not} work either. See Example \ref{Eg:1} for an illustration. Another possible attempt to extract a candidate is by deriving a discrete analogue of the Skorohod representation for $X^{(k)}$: one writes $X^{(k)}_t$ as the sum of a local martingale and a process of finite variation, then tries to show that the finite variational part converges in distribution to $\int_0^t\vec{n}\,(X_s)\,dL_s$. However, this has to be rigorously established. See Remark \ref{Rk:DiscretoLocalTime} (iv) below.

To the best of our knowledge, the question of discrete approximation to boundary local time of reflected diffusions has not even been rigorously addressed before. The main goal in this paper is to fill this gap. This paper is organized as follows:

In Section 2, we construct the discrete local time $L^{(k)}$ for RBM. This candidate is defined pathwise explicitly in \eqref{E:DiscreteLocalTime} (equivalently \eqref{E:DiscreteLocalTime1})  and is amenable to computer simulations. In Section 3, we state our main result, Theorem \ref{T:Main}, which is about weak convergence of joint laws $(X^{(k)},\,L^{(k)}) \to (X,\,L)$. Section 4 collects the key properties of transition density of $X^{(k)}$ including the local limit theorem,  Theorem \ref{T:LCLT_CTRW}, which is established in \cite{zqCwtF13a} with details in \cite{Fan14}. These properties will be used in the proof of Theorem \ref{T:Main} in Section 5. Extension of our main result to more general reflected diffusions is precisely stated in Theorem \ref{T:MainY} in Section 6.

\section{Discrete local time}

An important feature in our approach is that we incorporate geometric information of $\partial D$ in our approximation scheme. That is, besides approximating $D$ by $D^{(k)}$, we also approximate $\partial D$ by $\Lambda^{(k)}$, where for each $k \in\mathbb{N}$, $\Lambda^{(k)}$ is a partition of $\partial D$ into pieces of comparable sizes and diameters.
%This extra information is in a sense necessary for our explicit scheme, in view of Example \ref{Eg:1}. 
The choice of $\Lambda^{(k)}$ is specified by the following lemma.

\begin{lem}\label{L:DiscreteApprox_SurfaceMea}
Suppose $D$ is a bounded Lipschitz domain of $\R^d$. Then there exists a sequence of partitions $\{\Lambda^{(k)}\}_{k\in \mathbb{N}}$ of $\partial D$ and a constant $C\in (0,\infty)$ which depends only on $D$, such that (a), (b) and (c) below hold simultaneously:
\begin{enumerate}
\item[(a)]  %\begin{equation}\label{E:size_SurfaceMea}
            $ \sigma(\lambda)\leq C\,2^{-k(d-1)}$ for $\lambda\in \Lambda^{(k)}$ and $k \in\mathbb{N}$, where $\sigma$ is the surface measure on $\partial D$.
            %\end{equation}
\item[(b)]  %\begin{equation}\label{E:NumberOfPoints_SurfaceMea}
            $\sup_{x\in \bar{D}}\, \#\,\left\{\lambda\in \Lambda^{(k)}:\,\lambda \cap B(x,\,s)\neq \emptyset \right\} \leq C\,\left(2^ks\vee 1 \right)^{d-1}$ for $s\in(0,\infty)$ and $k \in\mathbb{N}$,
            %\end{equation}
            where $\# A $ is the cardinality of a finite set $A$ and $B(x,\,s)=\{y\in \R^d:\,|y-x|<s\}$.
\item[(c)] For any equi-continuous and uniformly bounded family $\F$ in $\mathcal{C}(\partial D)$, we have
            \begin{equation}\label{E:WeakConverge_SurfaceMea}
                \lim_{k\to \infty}\, \sup_{f\in\F} \sum_{\lambda\in \Lambda^{(k)}}\Big|\sup_{x\in \lambda}f(x)-\inf_{x\in \lambda}f(x)\Big|\,\sigma(\lambda)\,=0.
            \end{equation}
\end{enumerate}
\end{lem}

%\medskip
\noindent
The proof of Lemma \ref{L:DiscreteApprox_SurfaceMea} follows from an easy geometric argument which is basically a dyadic decomposition of $\partial D$. This proof can be found in \cite{zqCwtF13a}, in which a more general result about partitioning any rectifiable subsets of $\partial D$ is presented. \eqref{E:WeakConverge_SurfaceMea} implies that
$$\lim_{k\to\infty}\sum_{\lambda\in \Lambda^{(k)}}f(x_{\lambda})\,\sigma_{\lambda} = \int_{\partial D}f\,d\sigma$$ 
uniformly for $f\in\F$ and for all choices of $\{x_\lambda\}$ satisfying $x_{\lambda} \in \lambda$ for all $\lambda\in \Lambda^{(k)}$.

We are now describe our class of candidates for the desired discrete analogue for boundary local time, which is defined pathwise. % Such $z_{\lambda}$ exists since any point on $\partial D$ is of distance at most $\sqrt{1+M^2}\,2^{-k}$ to $D^{(k)}$.

\begin{definition}\rm\label{Def:DiscreteLocalTime}({\it Discrete local time})
Fix any $\alpha > \sqrt{1+M^2}$ where $M$ is the Lipschitz constant for $\partial D$. Associate each $\lambda\in \Lambda^{(k)}$ a non-empty subset $D^{(k)}_{\lambda}\subset D^{(k)}$ such that each $z\in D^{(k)}_{\lambda}$ is of distance at most $\alpha\,2^{-k}$ to $\lambda$.  Define, for each r.c.l.l. path $\omega:\,[0,\infty)\rightarrow D^{(k)}$ and $k\in \mathbb{N}$,
\begin{equation}\label{E:DiscreteLocalTime}
L^{(k)}_t(\omega) := \frac{1}{2}\int_0^t \sum_{\lambda\in \Lambda^{(k)}} \sum_{z\in D^{(k)}_{\lambda}}\dfrac{\1\{\omega(s)= z\}}{m_{k}(z)}\,\frac{\sigma(\lambda)}{\#D^{(k)}_{\lambda}}\;ds,
\end{equation}
where $m_{k}(x):=2^{-kd}\,v_{k}(x)/2d$ with $v_{k}(x)$ being the graph degree of the vertex $x\in D^{(k)}$. In particular, when $D^{(k)}_{\lambda}$ is a single point $\{z_{\lambda}\}$, then \eqref{E:DiscreteLocalTime} is reduced to
        \begin{equation}\label{E:DiscreteLocalTime_2}
            \frac{1}{2}\int_0^t \sum_{\lambda\in \Lambda^{(k)}}\dfrac{\1\{\omega(s)= z_{\lambda}\}}{m_{k}(z_{\lambda})}\,\sigma(\lambda)\,ds.
        \end{equation}
\end{definition}

\begin{remark}\label{Rk:DiscretoLocalTime}\rm
\begin{enumerate}
  \item[(i)] 
  Observe $D^{(k)}_{\lambda}$ is non-empty by the condition on $\alpha$, so that \eqref{E:DiscreteLocalTime} is well-defined. Note also that  $\#D^{(k)}_{\lambda}$ is abounded above by  some constant which depends only on the Lipschitz constant $M$. Furthermore, $\big\{D^{(k)}_{\lambda}:\,\lambda\in \Lambda^{(k)}\big\}$ can be flexibly chosen in such a way that $\partial^{(k)} :=\cup_{\lambda\in \Lambda^{(k)}}D^{(k)}_{\lambda}$ is equal to the graph boundary $\partial D^{(k)}$; in this case, $\#D^{(k)}_{\lambda}$ maybe larger than 1 for some $\lambda$ and we have to use \eqref{E:DiscreteLocalTime} rather than \eqref{E:DiscreteLocalTime_2}. 
 % \item[(ii)]  $\#D^{(k)}_{\lambda}\leq N$ for some constant $N$ which depends only on the Lipschitz constant $M$.
  \item[(ii)] Clearly, $L^{(k)}_t(\omega)$ is non-decreasing in $t$ and increases only when $\omega(t)\in \partial^{(k)}$. Hence
      \begin{equation*}
        L^{(k)}_t(\omega)= \int_0^t \1\{w(s)\in \partial^{(k)}\}\,dL^{(k)}_s(\omega).
      \end{equation*}
  \item[(iii)]  Intuitively, if the mass $\sigma(\lambda)$ of $\lambda$ is evenly distributed among elements in $D^{(k)}_{\lambda}$, then the total mass received by $z$ is given by
       $\sigma_k(z):=\sum_{\{\lambda:\,z\in D^{(k)}_{\lambda}\}}\sigma(\lambda)/\#D^{(k)}_{\lambda}$. The measure $\sigma_k$ on $\partial^{(k)}$ approximates $\sigma$ in the sense that
      $\lim_{k\to\infty}\sum_{z\in \partial^{(k)}}F(z)\,\sigma_k(z)\,=\,\int_{\partial D}F(z)\,\sigma(dz)$
      for any $F:\,D\to \R$ which is bounded and continuous on a neighborhood of $\partial D$. This is an immediate consequence of Lemma \ref{L:DiscreteApprox_SurfaceMea}. Moreover,  \eqref{E:DiscreteLocalTime} can be written as
        \begin{equation}\label{E:DiscreteLocalTime1}
            L^{(k)}_t(\omega) = \frac{1}{2}\int_0^t  \sum_{z\in \partial^{(k)}}\dfrac{\1\{\omega(s)= z\}}{m_{k}(z)}\,\sigma_k(z)\;ds.
        \end{equation}
  \item[(iv)]    In case $\partial^{(k)}$ is chosen to be $\partial D^{(k)}$,  which is always possible according to (i), then $X^{(k)}$ admits a pathwise decomposition analogous to \eqref{E:SkorokhodRep}:
      $$X^{(k)}_t=B^{(k)}_t+ \int_0^t \eta^{(k)}_s\,dL^{(k)}_s,$$
      where $B^{(k)}$ is the SRW
      % (continuous time or discrete time, according to $X^{(k)}$)
      on the whole lattice $ 2^{-k}\Z^d$, under the law of $X^{(k)}$; and $\eta^{(k)}$ is a $\F^{X^{(k)}}_t$-adapted process with values in $\R^d$. This ``Skorohod decomposition" can be used to study pathwise properties of $X^{(k)}$, but it will not play a role in our proof. 
      %We reserve discussions about its implications and the properties of $\eta^{(k)}$ in a future work. A related result can be found in  \cite{zqC93}, in which RBM is approximated by a sequence of RBMs on an increasing sequence of smooth domains.
\end{enumerate}
\end{remark}

\section{Main result and applications}

Recall that $X^{(k)}$ is the simple random walk on the graph $D^{(k)}$ moving at rate $d\,2^{2k}$, either continuous time or discrete time. In the latter case, time parameter is extended by interpolation as in \cite{BC08}. In each case, $X^{(k)}$ has stationary distribution $m_k$ stated in Definition \ref{Def:DiscreteLocalTime}. We denote by $\PP_{x_k}$ and $\PP_{m_k}$ the law of SRW $X^{(k)}$ starting from $x_k\in D^{(k)}$ and $m_k$ respectively. We also denote by $\P^{x}$ and $\P^{m}$ the law of RBM $X$ starting from $x\in \bar{D}$ and $m$ respectively, where $m$ is the uniform measure on $D$. %$\EE_{x_k}$, $\EE_{m_k}$, $\E^{x}$ and $\E^{m}$ denote the expectation with respect to $\PP_{x_k}$, $\PP_{m_k}$, $\P^{x}$ and $\P^{m}$ respectively. 
For a metric space $S$, we denote by $\mathcal{D}([0,T],S)$ the space of r.c.l.l. paths from $[0,T]$ to $S$ equipped with the Skorohod topology, and by $\mathcal{C}([0,T],S)$ the space of continuous paths equipped with uniform topology. Theorem \ref{T:Main} and Theorem \ref{T:MainY} are our main results.
% See also Theorem \ref{T:MainY} for a generalization to other reflected diffusions.

\medskip

\begin{thm}\label{T:Main}
Suppose $D$ is a bounded Lipschitz domain. Then for $T>0$, as $k\rightarrow \infty$ we have
\begin{enumerate}
\item[(i)]  $(X^{(k)},\,L^{(k)})$ under $\PP_{m_k}$ converges to $(X,\,L)$ in distribution in both $\mathcal{D}([0,T],\bar{D})\times \mathcal{C}([0,T],\R_+)$ and  $\mathcal{D}([0,T],\bar{D}\times \R_+)$, where $X$ is the reflected Brownian motion in $D$ with stationary initial distribution and $L$ is the boundary local time of $X$.
\item[(ii)]  If $x_k\in D^{(k)}$ converges to $x\in D$, then $(X^{(k)},\,L^{(k)})$ under $\PP_{x_k}$ converges to $(X,\,L)$ in distribution in both $\mathcal{D}([0,T],\bar{D})\times \mathcal{C}([0,T],\R_+)$ and  $\mathcal{D}([0,T],\bar{D}\times \R_+)$, where $X$ is the reflected Brownian motion in $D$ starting at $x$ and $L$ is the boundary local time of $X$.
\end{enumerate}
\end{thm}

\noindent
As an application, we consider the heat equation with general Robin boundary condition
    \begin{equation}\label{E:mixedBVP_drift}
        \left\{\begin{aligned}
        \dfrac{\partial u(t,x)}{\partial t} &= \frac{1}{2}\Delta u(t,x)   & &\qquad\text{on } (0,\infty)\times D  \\
        \dfrac{\partial u(t,x)}{\partial \vec{n}} &= g(t,x)\,u(t,x)+ h(t,x)  & &\qquad\text{on } (0,\infty)\times \partial D  %\\
        %u(0,x) &=\varphi(x) & &\qquad\text{on }  D,
        \end{aligned}\right.
    \end{equation}
and initial condition $f\in \mathcal{C}_b(D)$, where $g,\,h\in \mathcal{C}_b([0,\infty)\times \partial D)$  and $\mathcal{C}_b(E)$ denotes the space of bounded continuous functions on $E$. When $h=0$ this equation reduces to the classical Robin boundary problem. Using the Skorohod decomposition \eqref{E:SkorokhodRep} and It\^o formula, one obtains
a Feynman-Kac formula for the solution
\begin{equation}\label{F-Kac}
u(t,x)=\, \E^{x}\Big[f(X_t)\,e^{-\int^t_0g(t-s,X_s)\,dL_s}
- \int_0^t h(t-\theta,X_{\theta})\,e^{-\int^{\theta}_0 g(\theta-s,X_s)\,dL_s}\,dL_{\theta} 
\Big].
\end{equation}
See \cite[Proposition 2.17]{zqCwtF13a} for details of such a calculation.
Let  $G \text{ and } H \in  \mathcal{C}_b([0,\infty)\times \bar{D})$ be arbitrary continuous extensions of $g$ and $h$ respectively. 
Theorem \ref{T:Main} guarantees that 
\begin{equation*}\label{E:Application_Robin}
u_k(t,x_k):=\EE_{x_{k}}\Big[f(\omega(t))\,e^{-\int^t_0G(t-s,\omega(s))\,dL^{(k)}_s}
- \int_0^t H(t-\theta,\omega(\theta))\,e^{-\int^{\theta}_0 G(\theta-s,\omega(s))\,dL^{(k)}_s}\,dL^{(k)}_{\theta} 
\Big]
%u_k(t,x_k):= \EE_{x_{k}} \Big[\varphi(\omega(t)) \,\exp{\Big(-\int_0^tG(t-s,\,\omega(s))\,dL^{(k)}_s(\omega)\Big)} \Big]
\end{equation*}
converges to $u(t,x)$ whenever $x_k\to x\in \bar{D}$. Furthermore if $f \in C(\bar{D})$, then the convergence is uniform on $[a,b]\times \bar{D}$ for any compact interval $[a,b]\subset (0,\infty)$.

Since $L^{(k)}_s(\omega)$ increases only when $\omega(s)\in \partial^{(k)}:=\cup_{\lambda\in \Lambda^{(k)}}D^{(k)}_{\lambda}$, there is flexibility in the choice of $G$ and $H$.
%simply take $G(t,\,z):=g(t,\,z_{\lambda})$ in the following way: for $z\in \partial^{(k)}$, pick an arbitrary $\lambda$ such that $z\in D^{(k)}_{\lambda}$, then pick an arbitrary $z_{\lambda}$ in $\lambda$. The extension $H$ can be defined similarly. 
Hence Theorem \ref{T:Main} provides us with a convenient discrete approximation to the solution of \eqref{E:mixedBVP_drift}, using simple random walks and a decomposition of the boundary. 
Similar application of Theorem \ref{T:Main} also holds for elliptic equations (cf. \cite{vgP90}), using the probabilistic representation of the solutions.
%and lays down the foundation for Monte Carlo simulations for the solution to \eqref{E:mixedBVP_drift} using \eqref{E:Application_Robin}.

\medskip

The next two sections are devoted to the proof of Theorem \ref{T:Main}.
\vspace{-1em}
\section{Discrete heat kernel and local limit theorem}\label{LCLT}

In this section, we collect some fundamental properties of the transition density of random walks in domains.
%These properties are fundamental 
%that is needed in the proof of Theorem \ref{T:Main}. 
Most of these properties are proved in \cite{zqCwtF13a} for {\it biased random walks} which approximates RBM with drifts. See also \cite{Fan14} for detail of the calculations.
We consider  $D^{\eps}:= D\cap \eps \Z^d$ for $\epsilon>0$, and let $\partial D^{\eps}:= \{x\in D^{\eps}:\,v_{\eps}(x)<2d\}$ be the graph-boundary, where $v_{\eps}(x)$ is the degree of $x$ in $D^{\eps}$. We define $X^{\eps}$ to be the simple random walk (SRW) on $D^{\eps}$ moving at rate $d/\eps^{2}$, either continuous time or discrete time (as before, in the latter case, we extend time parameter by interpolation). Hence $X^{2^{-k}}$ in this section is the $X^{(k)}$ in Theorem \ref{T:Main}.

The transition density of $X^{\eps}$ with respect to measure $m_{\epsilon}(x):=\epsilon^{d}\,v_{\eps}(x)/2d$ is defined as
\begin{equation}\label{Def:pepsilon}
p^{\eps}(t,x,y) := \dfrac{\P^{x} (X^{\eps}_t=y)}{ m_\eps(y)}, \quad  t>0,\, x,\,y\in D^{\eps}.
\end{equation}
Clearly, $p^{\eps}$ is strictly positive and is symmetric in $x$ and $y$. It is proved in \cite{zqCwtF13a} that the transition density $p^{\eps}$ enjoys two-sided Gaussian bound and is jointly H\"older continuous
uniform in $\eps \in (0, \eps_0)$ for some $\eps_0>0$, and that $p^{\eps}$ converges to $p$ uniformly on compact subsets of $(0,\infty)\times \bar{D}\times \bar{D}$. In rigorous terms, we have the following four results. The important point is that the constants involved are uniform for $\epsilon$ small enough.

\begin{thm}\label{T:UpperHKE}(Gaussian upper bound)
There exist $C_k=C_k(d,D,T)\in (0,\infty)$, $k=1, 2$, and $\eps_0=\eps_0(d,D) \in (0, 1]$ such that for every $\eps\in(0,\eps_0)$ and
$x,y \in D^{\eps}$,
\begin{equation}\label{e:2.14}
p^{\eps}(t,x,y) \leq \dfrac{C_1}{(\eps\vee t^{1/2})^d}\,\exp\left(-\, C_2 \frac{|x-y|^2}{t}\right)
\quad \hbox{for }  t\in [\eps,T]\quad\text{and}
\end{equation}
\begin{equation}\label{e:2.15}
p^{\eps}(t,x,y) \leq \dfrac{C_1}{(\eps\vee t^{1/2})^d}\,\exp\left(-\,C_2 \frac{|x-y|}{ t^{1/2}}\right)
\quad \hbox{for } t\in (0, T].
\end{equation}
\end{thm}
%Observe that  \eqref{e:2.14} implies that \eqref{e:2.15} also holds for $t\in [\eps, T]$. 
%As an application, we have an estimate for the exit time for a ball by a standard argument (see \cite{BH91}) using the strong Markov property.

\begin{cor}\label{cor:ExitTime}(Exit time estimate)
There exist $C_k=C_k(d,D,T)\in (0,\infty)$, $k=1, 2$, and $\eps_0=\eps_0(d,D) \in (0, 1]$ such that for all $t\in(0,T]$,  $x\in D^{\epsilon}$, $\eta>0$ and $\epsilon\in(0,\epsilon_0)$,
\begin{equation}\label{E:ExitTime}
\P^{x}\Big(\sup_{s\leq t}|X^{\epsilon}_s-x| \geq \eta\Big) \leq C_1\,\exp\left(-\,\frac{C_2\,\eta}{(t^{1/2}\vee \epsilon)}\right).
\end{equation}
\end{cor}

\begin{thm}\label{T:LowerHKE}(Gaussian lower bound)
There exist $C_k=C_k(d,D,T)\in (0,\infty)$, $k=1, 2$, and $\eps_0=\eps_0(d,D) \in (0, 1]$ such that for every $\eps\in(0,\eps_0)$, $t\in (0,T]$ and $x,y \in D^{\eps}$,
\begin{equation}\label{e:2.19}
p^{\eps}(t,x,y) \geq \dfrac{C_1}{(\eps\vee t^{1/2})^d}\,\exp\left(- C_2 \frac{|x-y|^2}{t}\right).
\end{equation}
\end{thm}

\begin{thm}\label{T:HolderCts}(H\"older continuity)
There exist constants $\alpha (d,D,T)$, $\beta(d,D,T)$, $C(d,D,T)\in (0,\infty)$ and  $\eps_0(d,D)\in (0,1]$ such that for all $\eps\in(0,\eps_0)$, $(t,x,y),\,(t',x',y')\in (0,T]\times D^{\eps}\times D^{\eps}$, we have
\begin{equation}\label{E:HolderCts2}
|p^{\eps}(t,x,y)-p^{\eps}(t',x',y')| \leq C\, \dfrac{ (\,|t-t'|^{1/2}+ |x-x'|+ |y-y'|\,)^\alpha }
{(t\wedge t')^{(d+\beta)/2}}.
\end{equation}
\end{thm}

\begin{thm}\label{T:LCLT_CTRW}(Local limit theorem)
Let $p^{(k)}=p^{2^{-k}}$ be the transition density of $X^{(k)}$ with respect to $m_k$, and $p(t,x,y)$ be the transition density of the RBM with respect to Lebesque measure. Then we have
$$\lim_{k \to\infty} \sup_{t\in [a,b]} \sup_{x,y \in D^{(k)}}\Big|p^{(k)}(t,x,y)\,-\,p(t,x,y)\Big| =0 $$
for any compact interval $[a,b]\subset (0,\infty)$.
\end{thm}

The proofs for the above properties are standard once we establish a discrete analogue of a relative isoperimetric inequality in \cite[Theorem 5.5]{zqCwtF13a} for bounded Lipschitz domains. Details and stronger versions can be found in \cite{zqCwtF13a} and are omitted here. The following uniform estimate has a continuous analog. It is crucial to our proof of the main theorem.
\begin{lem}\label{L:q_near_I}
    There exist $C=C(d,D,T) \in (0,\infty)$ and $\eps_0=\eps_0(d,D) \in (0, 1]$ such that
    \begin{equation}
    \sup_{x\in D^{\eps}}\,\eps^{d-1}\,\sum_{y\in \partial D^{\eps}}p^{\eps}(t,x,y) \leq \frac{C}{\eps \vee t^{1/2}}
    \end{equation}
    for all $t\in (0,T]$ and $\eps\in(0,\eps_0)$.
\end{lem}

\begin{pf}
Fix $\theta\in (0,T]$. By the Gaussian upper bound in Theorem \ref{T:UpperHKE}, we have
\begin{eqnarray*}
&& \sum_{y\in \partial D^{\eps}}p^{\eps}(\theta,x,y) \\
&\leq& \frac{C_1}{(\eps \vee \theta^{1/2})^d}\,\sum_{y\in \partial D^{\eps}}\exp{\left(\frac{-|y-x|}{\eps \vee \theta^{1/2}}\right)} \\
&=& \frac{C_1}{(\eps \vee \theta^{1/2})^d}\int_0^{\infty}\#\,|\{y\in D^{\eps}:\,|f(y)|>r\}|\,dr \quad\text{where } f(y)=\1_{\partial D^{\eps}}(y)\,\exp{\left(\frac{-|y-x|}{\eps \vee \theta^{1/2}}\right)}\\
&=& \frac{C_1}{(\eps \vee \theta^{1/2})^d}\int_0^{1}\#\,|\{\partial D^{\eps}\cap B(x,\,(\eps \vee \theta^{1/2})(-\ln r))\}|\,dr \\%\quad (\text{since }f\leq 1)\\
&=& \frac{C_1}{(\eps \vee \theta^{1/2})^{d+1}}\int_0^{\infty}\#\,|\{\partial D^{\eps}\cap B(x,\,s)\}|\,\exp{\left(\frac{-s}{\eps \vee \theta^{1/2}}\right)}\,ds
\quad \text{where }s=(\eps\vee \theta^{1/2})(-\ln r),  \\
&\leq& \frac{C_1}{(\eps \vee \theta^{1/2})^d}\,\vee\,\frac{C_2}{\eps^{d-1}(\eps \vee \theta^{1/2})^{d+1}}\,\int_0^{\infty} s^{d-1}\,\exp{\left(\frac{-s}{\eps \vee \theta^{1/2}}\right)}\,ds \\
&\leq& \frac{1}{\eps^{d-1}}\,\left(\frac{C_1}{\eps \vee \theta^{1/2}}\,\vee \,\frac{C_2}{\eps \vee \theta^{1/2}}\,\int_0^{\infty} w^{d-1}e^{-w}dw \right) \quad
\text{where } w= \frac{s}{\eps \vee \theta^{1/2}}.
\end{eqnarray*}
Here $C_i$ are all constants which depend only on $d$, $D$ and $T$. Note that in the second last line, we used the fact, which follows from Lipschitz property of $\partial D$, that
$ \#\,|\{\partial D^{\eps}\cap B(x,\,s)\}| \leq C((s/\eps)^{d-1}\vee 1) $ for all $s>0$, for some $C=C(d,D)\in (0,\infty)$. The proof is now complete.
\end{pf}

\medskip

Recall $\partial^{(k)}$ in Remark \ref{Rk:DiscretoLocalTime}, which can be chosen to be $\partial D^{(k)}$. Lemma \ref{L:DiscreteApprox_SurfaceMea} implies that $ \#\,|\{\partial^{(k)} \cap B(x,\,s)\}| \leq C\,( 2^k\,s\vee 1)^{d-1} $ for some $C=C(d,D)\in (0,\infty)$. Hence the proof of Lemma \ref{L:q_near_I} implies

\begin{lem}\label{L:q_near_I2}
    There exist $C=C(d,D,T)\in (0,\infty)$ and $k_0=k_0(d,D)\in (0,\infty)$ such that
    \begin{equation}
    \sup_{x\in D^{(k)}}\,2^{k(d-1)}\,\sum_{y\in  \partial^{(k)}}p^{(k)}(t,x,y) \leq \frac{C}{2^{-k} \vee t^{1/2}}
    \end{equation}
    for all $t\in (0,T]$ and $k\geq k_0$, where $p^{(k)}$ is the transition density of $X^{(k)}$ with respect to $m_k$.
\end{lem}
This lemma is used crucially in the proof of Lemma \ref{L:Convergence}.

\section{Proof of main theorem}

In the following lemmas, we let $0\leq a\leq b$ and $\ell\in \mathbb{N}$ be arbitrary, and
$$\Delta_{\ell}[a,b] := \{(s_1,s_2,\cdots,s_{\ell}):\,a\leq s_1\leq s_2 \leq \cdots \leq s_{\ell} \leq b\}.$$
We also denote by $\mathcal{B}_b(E)$ the space of bounded measurable functions on $E$.

\begin{lem}\label{L:Identity1}
For $f\in \mathcal{B}_b(\partial D)$ and $x\in\bar{D}$, we have
\begin{eqnarray*}
\E^x\Big[\Big(\int_a^bf(X_s)\,dL_s\Big)^{\ell}\,\Big]&=&
\frac{{\ell}!}{2^{\ell}}\int_{\Delta_{\ell}[0,\,b-a]} \int_{\partial D}\cdots\int_{\partial D}\sigma(d y_1)\cdots\sigma(d y_{\ell})\,ds_1\cdots ds_{\ell}\\
&& \quad p(a+s_1,x,y_1)\,p(s_2,y_1,y_2)\cdots p(s_{\ell},y_{{\ell}-1},y_{\ell})\\
&& \quad f(y_1)\cdots f(y_{\ell})\;\rho(y_1)\cdots\rho(y_{\ell}).
\end{eqnarray*}
\end{lem}

\begin{pf}
Suppose $h\in \mathcal{B}_b( [0,T]\times\partial D)$. Then for $t\in[0,T]$ and  $x\in \bar{D}$, we have
\begin{equation}\label{E:expectation1}
\E^x\left[\int_0^th(s,X_s)dL_s\right] = \dfrac{1}{2}\int_0^t\int_{\partial D} h(s,y)\,p(s,x,y)\,\rho(y)\,\sigma(dy)\,ds.
\end{equation}
See \cite[Proposition 1.1]{vgP90} for the case when $D$ has $C^3$ boundary. For Libschitz boundary, the same proof goes through in view of  \cite{BH91}. The remaining computation is standard. We provide the detail here since it is used in the next lemma also.
By Fubini's Theorem and Markov property,
\begin{align*}
\E^x\Big[\Big(\int_0^t f(X_s)&\,dL_s\Big)^{\ell}\,\Big]= {\ell}!\,\E^x\int_{\Delta_{\ell}[0,t]}f(X_{s_{\ell}})\cdots f(X_{s_1})\,dL_{s_{\ell}}\cdots dL_{s_1}\\
&= {\ell}!\,\E^x\int^t_0\Big(\int_{\Delta_{{\ell}-1}[s_1,\,t]}f(X_{s_{\ell}})\cdots f(X_{s_2}) \, dL_{s_{\ell}}\cdots dL_{s_2}\Big)\,f(X_{s_1})\,dL_{s_1}\\
%&={\ell}!\,\E^x\int^t_0\Big(\int_{\Delta_{{\ell}-1}[0,\,t-s_1]}f(X_{s_{\ell}})\cdots f(X_{s_2})\,dL_{s_{\ell}}\cdots dL_{s_2}\Big)\circ\theta_{s_1}\;f(X_{s_1})\,dL_{s_1}\\
&={\ell}!\,\E^x\int^t_0\E^{x_{s_1}}\Big[\int_{\Delta_{{\ell}-1}[0,\,t-s_1]}f(X_{s_{\ell}})\cdots f(X_{s_2})\,dL_{s_{\ell}}\cdots dL_{s_2}\Big]\,f(X_{s_1})\,dL_{s_1}\\
&=\frac{{\ell}!}{2}\int^t_0 \int_{\partial D}p(s_1,x,y)\,g(y)\,\rho(y)\,\sigma(dy)\,ds_1 \quad\text{by }\eqref{E:expectation1},
\end{align*}
where 
%$\theta_s$ is the shift operator $(\theta_s(\omega))(t)=\omega(s+t)$, $\omega\in \mathcal{D}([0,\infty),\bar{D})$, and
%\begin{equation*}
$g(y)=\E^y\Big[\int_{\Delta_{{\ell}-1}[0,\,t-s_1]}f(X_{s_{\ell}})\cdots f(X_{s_2})\,dL_{s_{\ell}}\cdots dL_{s_2}\Big]\,f(y)$.
%\end{equation*}
By induction, the result for the case $a=0$ holds. The result also holds for $a>0$ by Markov property of $X$.
%$$\E^x\big[\big(\int_a^bf(X_s)\,dL_s\big)^{\ell}\,\big]=E^xE^{X_a}\big[\big(\int^{b-a}_0f(X_s)\,dL_s\big)^{\ell}\,\big] $$ by Markov property of the RBM $X$.
\end{pf}

By the same calculations and using the Makov property of $X^{(k)}_t$, we obtain
% =\EE_x\EE_{X^{(k)}_a}\Big[\Big(\int^{b-a}_0f(X^{(k)}_s)\,dL^{(k)}_s\Big)^{\ell}\Big]
\begin{lem}\label{L:Identity2}
For $f\in \mathcal{B}_b(D)$, $k \in \mathbb{N}$ and $x\in D^{(k)}$, we have
\begin{eqnarray*}
&&\EE_x\Big[\Big(\int_a^bf(X^{(k)}_s)\,dL^{(k)}_s\Big)^{\ell}\,\Big]
\\
&=&\frac{{\ell}!}{2^{\ell}}\int_{\Delta_{{\ell}}[0,\,b-a]}\sum_{\lambda_1\in\Lambda^{(k)}}\cdots\sum_{\lambda_{\ell}\in\Lambda^{(k)}}p^{(k)}(a+s_1,x,z_{\lambda_1}) p^{(k)}(s_2,z_{\lambda_1},z_{\lambda_2})\cdots p^{(k)}(s_{\ell},z_{\lambda_{{\ell}-1}},z_{\lambda_{\ell}})\\
&& \qquad \qquad \qquad \qquad \qquad \qquad \qquad f(z_{\lambda_1})\cdots f(z_{\lambda_{\ell}})\,\sigma(\lambda_1)\cdots \sigma(\lambda_{\ell})\,ds_1\cdots ds_{\ell}.
\end{eqnarray*}
\end{lem}

The next convergence result is the key in identifying subsequential limits of $(X^{(k)},L^{(k)})$.
\begin{lem}\label{L:Convergence}
For any $f\in \mathcal{B}_b(D)$ which is uniformly continuous in a neighborhood of $\partial D$,
\begin{equation}\label{E:Convergence1}
\lim_{k\to\infty}\EE_{x_k}\Big[\Big(\int^b_af(X^{(k)}_s)\,dL^{(k)}_s\Big)^{\ell}\,\Big]=\E^x\Big[\Big(\int^b_af(X_s)\,dL_s\Big)^{\ell}\,\Big]
\end{equation}
uniformly for $x\in\bar{D}$ and for any sequence $x_k\in D^{(k)}$ which converges to $x$. In particular,
\begin{equation}\label{E:Convergence2}
\lim_{k\to\infty}\EE_{m_k}\Big[\Big(\int^b_af(X^{(k)}_s)\,dL^{(k)}_s\Big)^{\ell}\,\Big]=\E^m\Big[\Big(\int^b_af(X_s)\,dL_s\Big)^{\ell}\,\Big].
\end{equation}
\end{lem}

\begin{pf}
It suffices to show the right hand side of the identities in Lemma \ref{L:Identity1} converges to that of Lemma \ref{L:Identity2} in the sense stated for \eqref{E:Convergence1}. We demonstrate the case $\ell=1$, as other cases can be proved in the same way. We want to show that
\begin{equation}\label{E:Convergence3}
\int_a^b\sum_{\lambda\in\Lambda^{(k)}}p^{(k)}(s,x_k,z_\lambda)\,f(z_\lambda)\,\sigma(\lambda)\,ds \to
\int_a^b\int_{\partial D}p(s,x,z)\,f(z)\,\sigma(dz)\,ds
\end{equation}
uniformly for $x\in\bar{D}$ and for any sequence $x_k\in D^{(k)}$ which converges to $x$. We first argue pointwise convergence. For fixed $s\in (a,b)$, the integrand (with respect to $ds$) converges by the local limit theorem (Theorem \ref{T:LCLT_CTRW}) and Lemma \ref{L:DiscreteApprox_SurfaceMea}. Hence by Lemma \ref{L:q_near_I2} and Lebesque dominated convergence theorem, we have \eqref{E:Convergence3} whenever $x_k\to x$.

By assumption on $f$, there exists $k_0$ large enough such that $f$ is uniformly continuous in a neighborhood of $\partial D$ which contains $\Lambda^{(k)}$ for all $k\geq k_0$. Besides, by interpolations (see, for example, \cite{zqCwtF13a}), $p^{(k)}$ can be viewed as an element in $\mathcal{C}([0,\infty)\times \bar{D}\times \bar{D})$. Now the desired uniform convergence follow from the pre-compactness of the sequence $\{g_k\} \subset \mathcal{C}(\bar{D})$, where
$g_k(x)=\int_a^b\sum_{\lambda\in\Lambda^{(k)}}p^{(k)}(s,x,z_\lambda)\,f(z_\lambda)\,\sigma(\lambda)\,ds$ is the left hand side of \eqref{E:Convergence3}. More precisely, uniform boundedness follows from Lemma \ref{L:q_near_I2}, while equicontinuity follows from the H\"older continuity of $p^{(k)}$ in Theorem \ref{T:HolderCts}.
\end{pf}

\bigskip

\begin{pf}\emph{of Theorem \ref{T:Main}:  }
By Lemma \ref{L:Identity2}, we have
\begin{eqnarray}\label{E:Main_Tight0}
\EE_x\Big[\Big(\int_a^bf(X^{(k)}_s)\,dL^{(k)}_s\Big)^{\ell}\,\Big]
&\leq& \frac{{\ell}!}{2^{\ell}}\,\| f\|^{\ell}\, C^{\ell} \int_{\Delta_{{\ell}}[0,b-a]}\frac{1}{\sqrt{(a+s_1)s_2\cdots s_{\ell}}}\,ds_1\cdots ds_{\ell} \notag\\
&\leq& \| f\|^{\ell} \, \frac{C^{\ell}\,{\ell}!}{\Gamma((\ell+2)/2)}(b-a)^{\ell/2}
\end{eqnarray}
for all $x\in D^{(k)}$ and $k\geq k_0=k_0(D) \in (0,\infty)$,  where $C=C(d,D,T)\in (0,\infty)$ and $\Gamma$ is the Gamma function. Taking $f\equiv1$, we obtain
\begin{equation}\label{E:Main_Tight}
\sup_{k\geq k_0}\sup_{x_k\in D^{(k)}}\EE_{x_k}\Big[\big|L^{(k)}_b-L^{(k)}_a\big|^{\ell}\Big]\leq C(b-a)^{\ell/2}
\end{equation}
for all $0\leq a\leq b\leq T$, where $k_0=k_0(D)\in (0,\infty)$ and $C=C(d,D,\ell,T)\in (0,\infty)$ are constants. By \eqref{E:Main_Tight} and the Kolmogorov-Centov tightness criteria (see \cite[Theorem 3.8.8]{EK86}), we obtain tightness of $\{L^{(k)}\}$ under $\{\PP_{x_k}\}$ in $\mathcal{C}([0,T],\R_+)$, where  $\{x_k\}$ is any sequence such that $x_k\in D^{(k)}$. Besides, \eqref{E:Main_Tight} clearly implies
\begin{equation}\label{E:Main_Tight2}
\sup_{k\geq k_0}\EE_{m_k}\Big[\big|L^{(k)}_b-L^{(k)}_a\big|^{\ell}\Big]\leq C(b-a)^{\ell/2}.
\end{equation}
Hence we also have the tightness of $\{L^{(k)}\}$ under $\{\PP_{m_k}\}$. By \cite[Lemma 2.1, Lemma 3.2]{BC08} and \cite[Remark 3.7]{BC11}, $\{X^{(k)}\}$ is tight in $\mathcal{D}([0,T],\bar{D})$ under both $\{\PP_{x_k}\}$ and $\{\PP_{m_k}\}$. The previous two sentences immediately imply tightness of $\{(X^{(k)},\,L^{(k)})\}$ in the product space $\mathcal{D}([0,T],\bar{D})\times \mathcal{C}([0,T],\R_+)$, under both $\{\PP_{x_k}\}$ and $\{\PP_{m_k}\}$. Tightness of $\{(X^{(k)},\,L^{(k)})\}$ in $\mathcal{D}([0,T],\bar{D}\times\R_+)$ also holds since the second component is continuous. It remains to identify subsequential limits.

We first consider subsequential limits in $\mathcal{D}([0,T],\bar{D}\times\R_+)$. Suppose, without loss of generality, that the full sequence $(X^{(k)},L^{(k)})$, under  $\{\PP_{m_k}\}$, converges in distribution to $(\tilde{X},\,\tilde{L})$ defined on some probability space $(\tilde{\Omega},\,\tilde{\F},\,\tilde{\P})$. Then results in \cite{BC08} implies that $\tilde{X}$ is the RBM under $\tilde{\P}$, because the map from $\mathcal{D}([0,T],\bar{D}\times\R_+)$ to $\mathcal{D}([0,T],\bar{D})$ which sends $(\omega_1,\omega_2)$ to $\omega_1$ is continuous (see problem 13 in \cite[Chapter 3]{EK86}). It remains to check that $\tilde{L}$ is the boundary local time of $\tilde{X}$ under $\tilde{\P}$.

We first show that $\tilde{L}$ is a PCAF of $\tilde{X}$. First, $\tilde{L}_t$ is continuous by \eqref{E:Main_Tight2}. This continuity then implies the convergence of finite dimensional distributions (see Theorem 7.8 in \cite[Chapter 3]{EK86})
$$(L^{(k)}_{t_1},\cdots,\,L^{(k)}_{t_m}) \to (\tilde{L}_{t_1},\cdots,\,\tilde{L}_{t_m}) \quad\text{as }k\to\infty$$
for all $0\leq t_1 < \cdots < t_m<\infty$. In particular, $\tilde{L}_0=0$ $\tilde{\P}$-a.s. By first considering rational numbers and then using continuity of $\tilde{L}$, we can check that $\tilde{L}_t$ is non-decreasing in $t$, since each of its prelimits is non-decreasing. Second, observe that $L^{(k)}$ is an additive functional by construction. Hence by convergence of joint distribution $(L^{(k)}_s,\,L_t^{(k)},\,L^{(k)}_s\circ \theta_t)$ for $t,\,s\geq 0$, we have $\tilde{L}_{t+s}(\omega)=\tilde{L}_{t}(\omega)+\tilde{L}_{s}(\theta_t\omega)$ a.s. for all $t,\,s\geq 0$. By continuity of $\tilde{L}$, we can strengthen the previous statement to obtain the additive property
$$\tilde{L}_{t+s}(\omega)=\tilde{L}_{t}(\omega)+\tilde{L}_{s}(\theta_t\omega), \quad t,\,s\geq 0,\,\tilde{P}\text{-a.s.}$$
Third,  $\tilde{L}_t$ is $\sigma(\tilde{X}_s:\,s\leq t)$ measurable by Skorohod representation theorem and the fact that  $L^{(k)}_t$ is $\sigma(X^{(k)}_s:\,s\leq t)$ measurable for all $k\in \mathbb{N}$ and $t\geq 0$.
These assert that $\tilde{L}$ is a PCAF of $\tilde{X}$.

Fix any $f\in \mathcal{C}_b(\partial D)$. Let $F\in \mathcal{C}_b(\bar{D})$ be any extension of $f$. The map $(\mu,\nu)\mapsto \int_0^{\cdot} F(\mu_s)d\nu_s$ is continuous from  $\mathcal{D}([0,T],\bar{D}\times\R_+)$ to  $\mathcal{D}([0,T],\R_+)$. Hence $\int_0^{\cdot}X^{(k)}_sdL^{(k)}_s \to \int_0^{\cdot}\tilde{X}_sd\tilde{L}_s$ in law in  $\mathcal{D}([0,T],\R_+)$. Since $\int_0^{t}\tilde{X}_sd\tilde{L}_s$ is continuous in $t$ by continuity of $\tilde{L}$, we have for all $t\geq 0$,
\begin{eqnarray*}
\tilde{E}\int_0^tf(\tilde{X}_s)\,d\tilde{L}_s
&=& \tilde{E}\int_0^tF(\tilde{X}_s)\,d\tilde{L}_s \\
&=& \lim_{k\to\infty} \EE_{m_k}\int_0^t F(X^{(k)}_s)\,dL^{(k)}_s \\
&=& \E_m\int_0^tf(X_s)\,dL_s \quad \text{by } \eqref{E:Convergence2} \\
&=& \frac{t}{2}\int_{\partial D} f(y)\,\sigma(dy) \quad \text{by } \eqref{E:expectation1}.
\end{eqnarray*}
By a standard monotone convergence argument, we have $\tilde{E}\int_0^tf(\tilde{X}_s)\,d\tilde{L}_s = \frac{t}{2}\int_{\partial D} f(y)\,\sigma(dy)$ for all $f\in \mathcal{B}_b(\partial D)$. Therefore, $\tilde{L}$ is the PCAF of $\tilde{X}$ associated with the measure $\sigma/2$ (see \cite[Appendix]{CF12}). By definition, $\tilde{L}$ is the boundary local time of $\tilde{X}$ under $\tilde{\P}$. The same arguments in the last three paragraphs work for subsequential limits of $(X^{(k)},L^{(k)})$ under $\{\PP_{x_k}\}$, using \eqref{E:Convergence1} rather than \eqref{E:Convergence2}. Therefore, sub-sequential limits in $\mathcal{D}([0,T],\bar{D}\times\R_+)$ are identified to be the same.
Finally, subsequential limits in $\mathcal{D}([0,T],\bar{D})\times \mathcal{C}([0,T],\R_+)$ can be identified in the same way. The proof is complete.
\end{pf}

\begin{example}\label{Eg:1} % Add figure
Let $D$ be the square with vertices $\{(1,0),\,(-1,0),\,(0,1),\,(0,-1)\}$ and  $C \in (\sqrt{2},\,3/\sqrt{2})$. Then $D^{C\,2^{-k}}\supset\partial D^{(k)}$ for all $k\in \mathbb{N}$ and for each $k$, the set $D^{C\,2^{-k}}\cap D^{(k)}$ remains the same for all such $C$. Arguing as in the proof of \eqref{E:Convergence3}, we have
\begin{eqnarray*}
\lim_{k\to \infty}\EE_{x_k}\big[A^{(k)}_{C\,2^{-k}}(t)\big]&=&\frac{3}{C\sqrt{2}}\,\E^x[L_t] \quad\text{and}\\
\lim_{k\to \infty}\EE_{x_k}\Big[\frac{1}{2\,(2^{-k})}\int_0^t \1\{X^{(k)}_s\in \partial D^{(k)}\}\,ds\Big]  &=& \frac{1}{\sqrt{2}}\, \E^x[L_t]
\end{eqnarray*}
whenever $x_k\to x$. Hence neither $A^{(k)}_{C\,2^{-k}}(t)$ nor $\frac{1}{2\,(2^{-k})}\int_0^t \1\{X^{(k)}_s\in \partial D^{(k)}\}\,ds$ is a suitable approximation to $L_t$.  It is clear that in the second case above, the factor $1/\sqrt{2}$ comes from the fact that only about $2^{k}$ points on each side of the square is used in the calculation of the left-hand side, while Definition \ref{Def:DiscreteLocalTime} asserts that about $2^{k}\sqrt{2}$ points should be used.
%A remedy using Definition \ref{Def:DiscreteLocalTime} can be achieved. For example,
\end{example}

\section{Extensions}

{\bf RBM with variable diffusion coefficient and gradient drift. }
%It is briefly mentioned in Section \ref{LCLT} that all estimates there hold for  {\it biased random walks} which approximates RBM with drifts. 
We now generalize our main result Theorem \ref{T:Main} to $\A$-reflected diffusions with
\begin{equation}\label{BMdrift}
\mathcal{A} : \,=\,\frac{a}{2}\,\Big(\Delta +\nabla h\cdot\nabla \Big)
\end{equation}
for some $a,\,h\in W^{1,2}(D)\cap\, \mathcal{C}(\bar{D})$ strictly positive. That is,  time-changed  Brownian motions with gradient drifts.
\eqref{BMdrift} corresponds to the general form $ \mathcal{A} := \frac{1}{2\,\rho}\,\nabla\cdot(\rho\,\textbf{a}\nabla)$ in \eqref{gen} with  $\textbf{a}(x)=a(x)\,I_{d\times d}$ and $\rho(x)=e^{2h(x)}/a(x)$.

To state the result precisely, we need to first construct a  {\it biased random walk} $Y^{\epsilon}$ on  $D^{\eps}= D\cap \eps \Z^d$. Define the symmetric weights (conductances) $\{\mu_{xy}:\,x,y\in D^{\epsilon} \text{ adjacent}\}$ by two steps: First, assign for every $x\in D^{\epsilon}\setminus \partial D^{\epsilon}$ and $i=1,2,\cdots,d$,
\begin{eqnarray*}
	 \mu_{x,x+\epsilon \vec{e_i}} &:=& (1+h(x+\epsilon \vec{e_i})-h(x))\,\left(\frac{e^{2 h(x)}+e^{2 h(x+\epsilon \vec{e_i})}}{2}\right)\,\frac{\epsilon^{d-2}}{2}\\
 \mu_{x,x-\epsilon \vec{e_i}} &:=& (1+h(x)-h(x-\epsilon \vec{e_i}))\,\left(\frac{e^{2 h(x)}+e^{2 h(x-\epsilon \vec{e_i})}}{2}\right)\,\frac{\epsilon^{d-2}}{2},
\end{eqnarray*}
so that $\mu_{xy}=\mu_{yx}$ for all $x,y\in D^{\epsilon}\setminus \partial D^{\epsilon}$. Second, extend to define
$$\mu_{xy} \triangleq
\begin{cases}
\mu_{yx}, &\text{ if  } x\in \partial D^{\epsilon} \text{ and }\,y\in D^{\epsilon}\setminus \partial D^{\epsilon} \text{ are adjacent}\\
\epsilon^{d-2}/2, &\text{ if  } x,y\in \partial D^{\epsilon}  \text{ are adjacent}.
\end{cases}$$
Now $\mu_{xy}=\mu_{yx}$ for all $x,y\in D^{\epsilon}$. Let $\mu_{\epsilon}(x):= \sum_{y }\mu_{xy}$. 

%Consider the bilinear form $(\D^{(\epsilon)},\, l^2(m_{\epsilon}))$ defined by
%\begin{eqnarray*}
%	\D^{(\epsilon)} (f,g)  &\triangleq &  \dfrac{1}{2} \sum_{x,y\in D^{\epsilon}}(f(y)-f(x))(g(y)-g(x))\,\mu_{xy},
%\end{eqnarray*}
%Clearly, $(\D^{(\epsilon)},\, l^2(m_{\epsilon}))$ is a regular Dirichlet form in $l^2(m_{\epsilon})$. By Theorem 7.2.1 of \cite{FOT94}, there is a unique associated $m_{\epsilon}$-symmetric Markov process $X^{\epsilon}$.

\begin{definition}\rm\label{Def:Yeps}
Let $Y^{\epsilon}$ be the {\it biased random walk} on $D^{\epsilon}$ with jump rate  $\lambda_{\epsilon}(x) = a(x)d/\epsilon^{2}$ and one step transition probabilities $p_{xy} \triangleq \mu_{xy}/ \mu(x)$.	As before,  $Y^{\epsilon}$ can be either continuous time or discrete time. In the latter case, we extend time parameter by interpolation. 
We also let $Y$ be the reflected diffusion with generator given by \eqref{BMdrift}.
It is easy to check that $Y^{\epsilon}$ and $Y$ are symmetric with respect to $m_{\epsilon}(x):= \mu(x)/\lambda_{\epsilon}(x)$ and  $m(x):=e^{2h(x)}/a(x)$ respectively. 
%We call $Y^{\epsilon}$ the $\epsilon$-approximation of the reflected diffusion $Y$. 
\end{definition}

Our generalization to Theorem \ref{T:Main} is precisely stated below. It is remarkable that {\em the same} $L^{(k)}$ in Definition \ref{Def:DiscreteLocalTime} can be used. As before, $Y^{(k)}=Y^{2^{-k}}$ and $m_k=m_{2^{-k}}$.

\begin{thm}\label{T:MainY}
	Suppose $D\subset \R^d$ is a bounded Lipschitz domain. Suppose $a,\,h\in W^{1,2}(D)\cap\, \mathcal{C}(\bar{D})$ are strictly positive.  Let $\PP_{x_k}$ and $\PP_{m_k}$ be the laws of $Y^{(k)}$ starting from $x_k\in D^{(k)}$ and  $m_{k}(x)$ respectively. Let $\P^{x}$ and $\P^{m}$ be the laws of $Y$ starting from $x\in \bar{D}$ and $m(x):=e^{2h(x)}/a(x)$ respectively. For every $T>0$, as $k\rightarrow \infty$, the followings hold:
	\begin{enumerate}
		\item[(i)]  $(Y^{(k)},\,L^{(k)})$ under $\PP_{m_k}$ converges to $(Y,\,L)$ in distribution in both $\mathcal{D}([0,T],\bar{D})\times \mathcal{C}([0,T],\R_+)$ and  $\mathcal{D}([0,T],\bar{D}\times \R_+)$, where $Y$ has stationary initial distribution $m(x)dx$ and $L$ is the boundary local time of $Y$.
		\item[(ii)]  If $x_k\in D^{(k)}$ converges to $x\in D$, then $(Y^{(k)},\,L^{(k)})$ under $\PP_{x_k}$ converges to $(Y,\,L)$ in distribution in both $\mathcal{D}([0,T],\bar{D})\times \mathcal{C}([0,T],\R_+)$ and  $\mathcal{D}([0,T],\bar{D}\times \R_+)$, where $Y$ starts at $x$ and $L$ is the boundary local time of $Y$.
	\end{enumerate}
\end{thm}

\begin{pf}
Suppose $a,\,h\in W^{1,2}(D)\cap\, \mathcal{C}(\bar{D})$ strictly positive. Then from Theorem 2.2.20 in \cite{Fan14}, $Y^{\epsilon}$ converges weakly to  $Y$. Moreover, let $q^{\epsilon}(t,x,y)$ be the transition density of $Y^{\epsilon}$ with respect to  $m_{\epsilon}(x)$. Then $q^{\epsilon}(t,x,y)$ converges locally uniformly to the transition density of $Y$ with respect to $m(x)$. In other words, the local central limit theorem holds. Furthermore, all estimates in Section \ref{LCLT} hold for $q^{\epsilon}(t,x,y)$ (see Section 2.2.5 in \cite{Fan14}). Now by the same argument used to prove Theorem \ref{T:Main}, it is straightforward to check that Theorem \ref{T:Main} remains true even if we generalize from RBM to reflected diffusions with generator \eqref{BMdrift}.
\end{pf}

\medskip

\begin{remark}\rm
Nearest neighbor random walk
approximations, such as the SRW $X^{\epsilon}$ and the biased random walk $Y^{\epsilon}$ in definition \ref{Def:Yeps}, are very desirable from the point of view of computer simulation and numerical algorithm. 
Nonetheless, it require a nontrivial amount of extra work to generalize Theorem \ref{T:Main} or Theorem \ref{T:MainY} 
to general reflected diffusions (such as when the matrix $\textbf{a}(x)$ is not of diagonal form). It seems, in view of results in  \cite{SZ97, BK08}, that nearest neighbor random walk approximations becomes highly nontrivial even for symmetric diffusions on $\R^d$. The Markov chain approximations in \cite{SZ97, BK08}  are not nearest neighbor. One can expect that, due to regularity issues on the boundary, nearest neighbor approximations of general reflected diffusions are more challenging to establish.
\end{remark}

\medskip

\noindent
{\bf Other extensions. }The idea in this paper can be easily extended to construct discrete approximations to other  positive continuous additive functionals (PCAF), such as the local time on any ($d-1$)-dimensional rectifiable subset in $\bar{D}$, such as an open subset of $\partial D$, the slit $[0,1)\times\{0\}$ in the unit disc, etc. 
The sequence $2^{-k}$ for the lattice size in this paper is chosen to follow that in \cite{BC08}.  Generalization of results in \cite{BC08} and this paper to any sequence which tends to zero is left to the readers. The fact that all estimates in Section 4 hold for  $\epsilon>0$ small enough will be useful.

\medskip

\begin{flushleft}
{\bf Acknowledgements}

The author thanks Amarjit Budhiraja (University of North Carolina), Krzysztof Burdzy and Zhen-Qing Chen (University of Washington) for thoughtful remarks.  Discussions with On Shun Pak (Princeton University) about scientific applications of reflected diffusions are appreciated. This research is partially supported by Army Research Office W911NF-10-1-0158.	
\end{flushleft}

\footnotesize{{\sc
		\bigskip

		\noindent
		W-T. Fan\\
		Department of Mathematics\\
		University of Wisconsin\\
		Madison, WI 53706, USA\\
		www.math.wisc.edu/~louisfan\\
		email: louisfan@math.wisc.edu

	}}
\end{document}